\documentclass[12pt,a4paper,reqno]{amsart}

\usepackage{amssymb}
\usepackage{amscd}
\usepackage[pdftex,pdfpagelabels]{hyperref}
\usepackage{enumerate}
\usepackage{comment}
\usepackage{graphicx}
\usepackage{siunitx}
\usepackage{tikz-cd}
\usepackage{stix}
\usepackage{bm}
\DeclareMathAlphabet\mathbfcal{LS2}{stixcal}{b}{n}
\numberwithin{equation}{section}
\usepackage{cleveref}
\usepackage[pagewise]{lineno}

\usepackage{mathtools}
\usepackage[tableposition=top]{caption}
\usepackage{booktabs,dcolumn}

\DeclareFontFamily{OT1}{rsfs}{}
\DeclareFontShape{OT1}{rsfs}{n}{it}{<-> rsfs10}{}
\DeclareMathAlphabet{\mathscr}{OT1}{rsfs}{n}{it}

\theoremstyle{plain}

\newtheorem{theorem}{Theorem}[section]
\newtheorem{proposition}[theorem]{Proposition}
\newtheorem{lemma}[theorem]{Lemma}
\newtheorem{corollary}[theorem]{Corollary}

\theoremstyle{definition}

\newtheorem{definition}[theorem]{Definition}

\newtheorem{example}[theorem]{Example}

\renewcommand\P{\mathbb{P}}
\newcommand\E{\mathbb{E}}
\newcommand\I{\mathbb{I}}

\newcommand\R{\mathbb{R}}
\newcommand\Z{\mathbb{Z}}

\renewcommand\H{\mathbb{H}}
\newcommand\h{\mathbf{h}}

\newcommand\Q{\mathbb{Q}}

\newcommand\SD{{\operatorname{SD}}}
\newcommand\eps{\varepsilon}

\begin{document}

\title[Sum-difference exponents and rational complexity]{Sum-difference exponents for boundedly many slopes, and rational complexity}

\author{Terence Tao}
\address{UCLA Department of Mathematics, Los Angeles, CA 90095-1555.}
\email{tao@math.ucla.edu}



\begin{abstract}  The dimension of Kakeya sets can be bounded using sum-difference exponents $\SD(R;s)$ for various sets of rational slopes $R$ and output slope $s$; the arithmetic Kakeya conjecture, which implies the Kakeya conjecture in all dimensions, asserts that the infimum of such exponents is $1$.  The best upper bound on this infimum currently is $1.67513\dots$.  In this note, inspired by numerical explorations from the tool \texttt{AlphaEvolve}, we study the regime where the cardinality of the set of slopes $R$ is bounded.  In this regime, we establish that these exponents converge to $2$ at a rate controlled by the \emph{rational complexity} of $s$ relative to $R$, which measures how efficiently $s$ can be expressed as a rational combination of slopes in $R$.
\end{abstract}  

\maketitle

\section{Introduction}

\subsection{The arithmetic Kakeya conjecture}

Define a \emph{slope} to be an element $r$ of the projective rational\footnote{One could also formulate the arithmetic Kakeya conjecture in other fields than the rationals, and in fact most of the results here apply to arbitrary infinite fields, with minor modifications in the case where the field has positive characteristic.  But in this paper we shall restrict attention for simplicity to the rational case, which is the case of interest for applications to the Kakeya problem.} line $\Q \cup \{\infty\}$.  We then define the projection operators $\pi_r \colon \Q \times \Q \to \Q$ by setting
$$ \pi_r(x,y) \coloneqq x + ry$$
for $r \neq \{\infty\}$, and
$$ \pi_\infty(x,y) \coloneqq y.$$
Given a finite set $R$ of slopes and a further slope $s$ not in $R$, we define the \emph{sum-difference constant} $\SD(R;s)$ to be the least exponent such that the bound
\begin{equation}\label{sdef}
\H[ \pi_s(X,Y) ] \leq \SD(R;s) \max_{r \in R} \H[\pi_r(X,Y)]
\end{equation}
holds for $\Q$-valued random variables\footnote{In this paper, all random variables are understood to be discrete, and in fact take only finitely many values.} $X,Y$ (not necessarily independent), where
$$ \H[X] \coloneqq \sum_{x} \h(\P(X=x))$$
is the Shannon entropy of a random variable $X$, with $\h(t) \coloneqq t \log \frac{1}{t}$ using the convention $\h(0) \coloneqq 0$.  Thus for instance
$$\H[ X-Y ] \leq \SD(\{0,1,\infty\};-1) \max\left( \H[X], \H[X+Y], \H[Y] \right)$$
whenever $X,Y$ are (possibly dependent) $\Q$-valued random variables.

The quantity $\SD(R;s)$ can equivalently be defined as the least exponent such that the bound
\begin{equation}\label{pise}
|\pi_{s}(E)| \leq (\max_{r \in R} |\pi_s(E)|)^{\SD(R;s)}
\end{equation}
for all finite non-empty $E \subset \Q$, where $|A|$ denotes the cardinality of a finite set $A$; see \cite{green}.  However, it will be convenient in this paper to work with the entropy formulation, in order to take advantage of the ``entropic Pl\"unnecke--Ruzsa calculus'' that are founded on the Shannon entropy inequalities.

It is easy to see that one has the projective invariance
$$ \SD(\phi(R); \phi(s)) = \SD(R; s)$$
for any projective transformation $\phi \colon \Q \cup \{\infty\} \to \Q \cup \{\infty\}$, that is to say a map of the form $\phi(r) \coloneqq \frac{ar+b}{cr+d}$ for some $a,b,c,d \in \Q$ with $ad-bc \neq 0$, with the usual conventions when $r$ is infinite or $cr+d$ vanishes.  For instance, by using a dilation transformation, we have
$$ \SD(\{0,1,\infty\}; s) = \SD(\{0,-1/s,\infty\}; -1)$$
for any slope $s$ other than $0,1,\infty$.
In the literature it is conventional to use this $3$-transitive projective symmetry to normalize $s = -1$ and $0, \infty \in R$ (assuming that $|R|\geq 2$ of course), though in this paper it will be more convenient to adopt the normalization $0,1,\infty \in R$ (assuming $|R| \geq 3$).

It is easy to see that $\SD(R;s)=\infty$ when $|R|<2$.
From the entropy inequality $\H[X-Y] \leq \H[X,Y] \leq \H[X]+\H[Y]$ we see that $\SD(R;-1) \leq 2$ when $0,\infty \in S$, and it is easy to see (using the uniform distribution on a long arithmetic progression) that we have equality when $R = \{0,\infty\}$.  By projective invariance, this implies that $\SD(R;s)=2$ whenever $|R|=2$.  As these quantities are clearly non-decreasing in $R$, we then have the trivial bound
\begin{equation}\label{sdf-upper}
\SD(R;s) \leq 2
\end{equation}
for $|R| \geq 2$.  

Improvements upon \eqref{sdf-upper} directly lead to improved upper bounds on the dimension of Kakeya and Nikodym sets in high dimensions.  Indeed, in \cite{bourgain} it was observed that a bound of the form $\SD(R;-1) \leq \alpha$ implies that Kakeya and Nikodym sets in dimension $d$ have (upper) Minkowski dimension\footnote{It was recently observed by Thomas Bloom (private communication) that, by combining Bourgain's arguments with the recent quantitative progress on Szemer\'edi's theorem by Leng, Sawhney, and Sah \cite{lss}, that one also obtains this bound for the Hausdorff dimension as well.  For the corresponding results for packing dimension, see \cite{cowen}.} at most $\frac{d-1}{\alpha}+1$.  In particular, if one can establish the \emph{arithmetic Kakeya conjecture}
\begin{equation}\label{infr}
\inf_R \SD(R;-1) = 1
\end{equation}
then this would imply that Kakeya and Nikodym sets in $\R^d$ have full Minkowski and Hausdorff dimension for all $d$.  This is currently only known for $d \leq 3$ \cite{wang}. We refer the reader to \cite{green}, \cite{cowen}, \cite{pohoata} for a discussion of this conjecture (and other equivalent forms of it), and its connection with other variants of the Kakeya conjecture.  

Nontrivial progress towards the arithmetic Kakeya conjecture was first obtained in  \cite{bourgain}, who in our notation showed that 
$$ \SD(\{0,1,\infty\};-1)  \leq 2 - \frac{1}{13} = 1.923\dots,$$
and used this to obtain new bounds on the Kakeya conjecture in high dimensions.
Further improvements were then obtained in \cite{katz-tao}, \cite{katz-tao-new}. For instance, it is known that
$$ 1.77898 \leq \SD(\{0,1,\infty\};-1) \leq 2 - \frac{1}{6} = 1.833\dots$$
and
$$ 1.668 \leq \SD(\{0,1,2,\infty\};-1) \leq 2 - \frac{1}{4} = 1.75,$$
with the upper bounds established in \cite{katz-tao}, and the lower bounds in \cite{lemm}, \cite{gdm} respectively.
At present, the best upper bound known towards \eqref{infr} is
$$\inf_R \SD(R;-1) \leq 1.67513\dots;$$
see \cite{katz-tao-new}.  

\subsection{Asymptotic behavior}

Informally, the arithmetic Kakeya conjecture asserts that in the asymptotic regime where the number of slopes $R$ is large, the constants $\SD(R;s)$ converge to $1$.  Here we consider a complementary regime, in which the number of slopes $R$ is fixed, but we instead let the elements of $R$ (or $s$) vary.  Our main results assert, roughly speaking, that the behavior of these constants is determined by the \emph{rational complexity} of $s$ relative to $R$, with the constants approaching $2$ as it becomes harder to express $s$ in terms of a rational expression of the $R$.  We give a (slightly artificial) definition of this quantity, restricting attention to the normalized setting $\{0,1,\infty\} \subset R$ for simplicity.

\begin{definition}[Rational complexity] Given a family of slopes $R = \{0,1,\infty,r_1,\dots,r_k\}$ containing $0,1,\infty$ and a further slope $s$ not in $R$, we define the \emph{rational complexity} $D = D(R;s)$ of $s$ relative to $R$ to be the least natural number $D$ for which one has a representation of the form
\begin{equation}\label{rat}
 s = \frac{P(r_1,\dots,r_k)}{Q(r_1,\dots,r_k)}
 \end{equation}
where $P,Q$ are polynomials of degree at most $D$ with integer coefficients of magnitude at most $2^D$, with $Q(r_1,\dots,r_k)$ non-zero; this complexity is finite since $s$ is rational.
\end{definition}

Informally, if the complexity of $s$ with respect to $R$ is equal to $D$, then $s$ can be expressed in terms of the slopes in $R$ by a rational expression whose length (when expressed as string of characters) is comparable to $D$.  The rational complexity is reminiscent of the \emph{arithmetic circuit complexity} of $s$ in terms of $R$, but with the key difference that the circuit must take the specific rational form \eqref{rat}.

\begin{example}  $D(\{0,1,\infty\};s)$ is the least natural number for which one can express $s$ as $a/b$ where $a,b$ are integers of magnitude at most $2^D$.  In particular, if $a,b$ are coprime then
$$D \left( \{0,1,\infty\};\frac{a}{b} \right) \asymp \log(2+|a|+|b|).$$
\end{example}

Our main result, proven in Section \ref{main-sec}, is then as follows.  

\begin{theorem}\label{main-2}  Let $R$ be a finite family of slopes containing $0,1,\infty$ of some cardinality $k+3$, and let $s$ be a slope not lying in $R$.
\begin{itemize}
    \item[(i)] (Three slopes) If $k=0$, then
\begin{equation}\label{2ab} 2 - \frac{c_2}{D(R;s)} \leq \SD\left(R; s\right) \leq 2 - \frac{c_1}{D(R;s)} 
\end{equation}
    for some absolute constants $c_2 > c_1 > 0$.
    \item[(ii)]  (Many slopes) In general, one has
    \begin{equation}\label{sam} 2 - \frac{C_k \log(2+D(R;s))}{D(R;s)} \leq \SD(R; s) \leq 2 - \frac{c_k}{D(R;s)^{k+1}} 
    \end{equation}
    for some absolute constants $c_k, C_k > 0$.
\end{itemize}
\end{theorem}

\begin{figure}
    \centering
\centerline{\includegraphics[width=\linewidth]{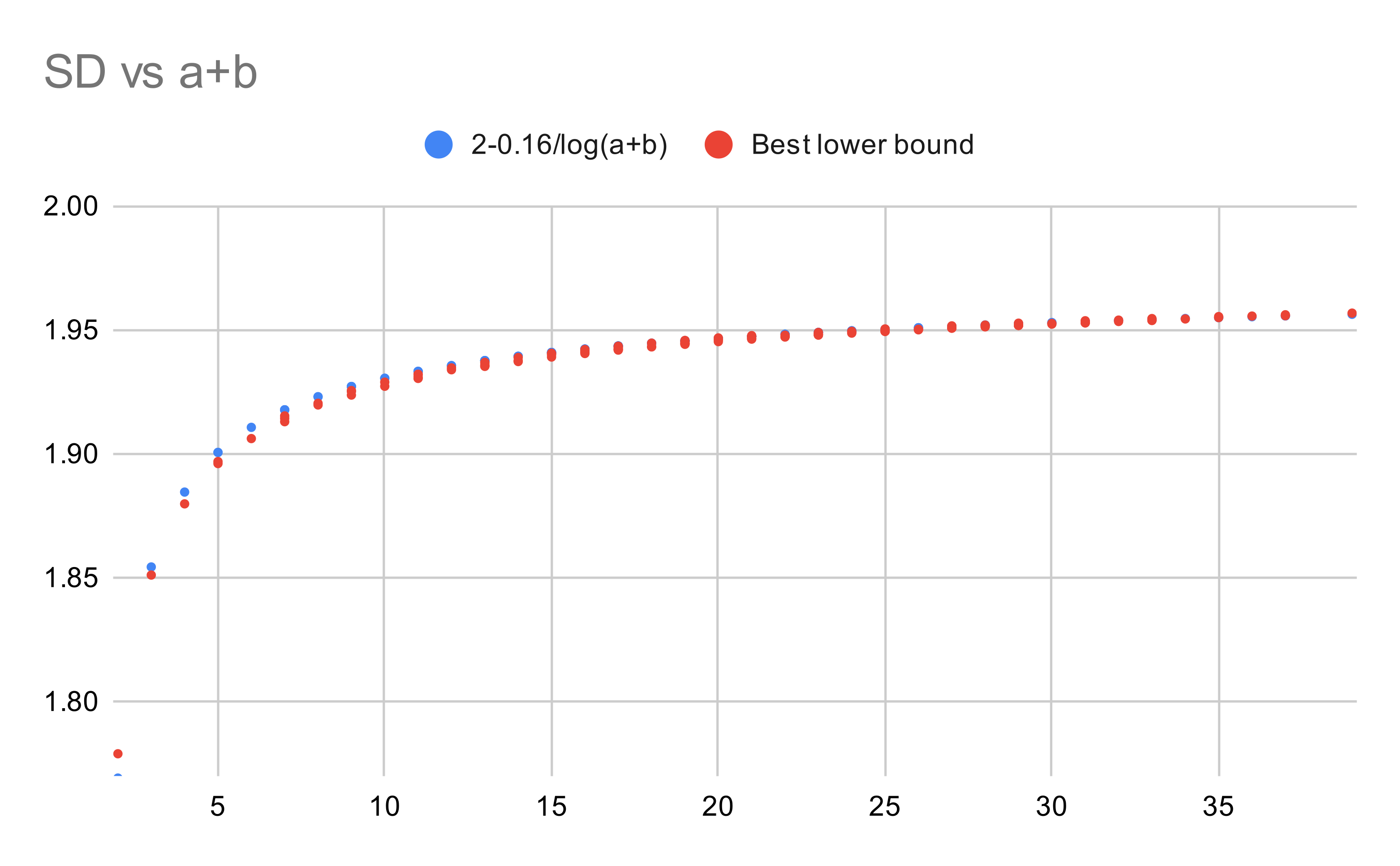}}
    \caption{Lower bounds on $\SD(\{0,1,\infty\}; s)$ obtained by \texttt{AlphaEvolve} for various $s=a/b$, plotted against $2 - 0.16 / \log(|a|+b|)$.  The horizontal axis is $|a|+|b|$ (so in some cases, multiple fractions $a/b$ are plotted on a single vertical line).}
    \label{fig:kak}
\end{figure}

The logarithmic convergence in \eqref{2ab} was suggested to us by experiments in \cite{gdm} using \texttt{AlphaEvolve} to obtain lower bounds for $\SD(\{0,1,\infty\}; s)$ for various slopes $s$.  This data was of low accuracy, as \texttt{AlphaEvolve} could only provide lower bounds and not upper bounds for these quantities; nevertheless, a logarithmic decay was numerically evident (see Figure \ref{fig:kak}), and furthermore the approximate discrete gaussian shape of the joint distributions of the random variables $X,Y$ obtained by this tool (see \cite[Figure 18]{gdm}) suggested an approach to make the lower bound in \eqref{2ab} rigorous.  Once this was accomplished, the author was also able to obtain the matching upper bound in \eqref{2ab} using the entropic Pl\"unnecke--Ruzsa calculus.   A modification of these arguments, with some inefficiencies, then gives \eqref{sam}. We tentatively conjecture that the bounds in \eqref{sam} can be improved to be of the form \eqref{2ab} for all $k$, not just $k=0$ (possibly after some slight adjustments to the definition of rational complexity).

In the limit as $a/b$ converges to some real number $\alpha$ and $b \to \infty$, we can obtain a more precise lower bound (also suggested by the aforementioned \texttt{AlphaEvolve} numerics) as follows.

\begin{theorem}[Continuous limit]\label{cts-limit}  Let $\alpha$ be a real number, and let $a,b$ be coprime integer parameters with $b \to \infty$ and $a/b \to \alpha$.  Then
\begin{equation}\label{up}
 \SD\left(\{0,1,\infty\};\frac{a}{b}\right) \geq 2 - \frac{c_\alpha+o(1)}{\log b}
 \end{equation}
 as $b\to \infty$ and $a/b \to \alpha$, where $0 < c_\alpha < \infty$ is the infimum of all quantities
\begin{equation}\label{hf}
 h_2(f_{(\alpha)}) - \max(h(f_0), h(f_\infty), h(f_\alpha))
 \end{equation}
where $f \colon \R^2 \to \R^+$ ranges over smooth compactly supported functions with total mass $\int_{\R^2} f(x,y)\ dx dy$ equal to $1$, $f_{(\alpha)} \colon \R \times [0,1) \to \R^+$ is a ``folded'' version
\begin{equation}\label{falpha-def}
 f_{(\alpha)}(x,y) \coloneqq \sum_{j \in \Z} f(x-\alpha j, y+j)
 \end{equation}
of $f$, $f_0, f_\infty, f_1 \colon \R \to \R^+$ are the projections
\begin{align}
f_0(x) &\coloneqq \int_\R f(x,y)\ dy\label{f0-def}\\
f_\infty(y) &\coloneqq \int_\R f(x,y)\ dx\label{finfty-def}\\
f_1(z) &\coloneqq \int_\R f(x,z- x)\ dx,\label{f1-def}
\end{align}
the two-dimensional differential entropy $h_2(f_{(\alpha)})$  of $f_{(\alpha)}$ is defined as
$$h_2(f_{(\alpha)}) \coloneqq \int_0^1 \int_\R \h(f_{(\alpha)}(x,y))\ dx dy,$$
and the one-dimensional differential entropies $h(f_r)$ for $r=0,1,\infty$ is defined as
$$ h(f_r) \coloneqq \int_\R \h(f_r(x))\ dx.$$
\end{theorem}

We establish this result in \Cref{cts-sec}.
We tentatively conjecture that the upper bound in \eqref{up} is in fact an asymptotic equality, so that the asymptotic behavior of
$ \SD\left(\{0,1,\infty\}, \frac{a}{b}\right)$ is controlled not only by the rational complexity (as represented by the $\log b$ denominator), but also by the variational quantity $c_\alpha$ appearing in the numerator.  It is not clear whether this quantity $c_\alpha$ can be computed exactly; numerically, two-dimensional gaussians are reasonably good candidates for $f$, but in practice they do not extremize the functional \eqref{hf} precisely.

\subsection{Notation}

We use the asymptotic notation $X = O(Y)$, $X \ll Y$, or $Y \gg X$ to denote the assertion that $|X| \leq CY$ for some absolute constant $C$; if we need this implied constant $C$ to depend on some fixed quantities (such as the number $k$ of slopes), we will indicate this in the text.

When there is no possibility of ambiguity, we omit parentheses from pairs $(X,Y)$ of random variables, for instance abbreviating $\H[(X,Y)]$ as $\H[X,Y]$.  Given a random variable $Y$ taking values in some set $S$ and some function $f \colon S \to \R$, we define the expectation
$$ \E_{Y=y} f(y) \coloneqq \sum_y \P(Y=y) f(y)$$
where $y$ ranges over the essential range of $Y$.  While this expression could also be abbreviated as $\E f(Y)$, it will be notationally useful to distinguish between the random variable $Y$ and the possible values $y$ that this variable could take.  For instance, with this notation, the \emph{conditional entropy} $\H[X|Y]$ of one random variable $X$ with respect to another $Y$ can now be defined by the formula
$$ \H[X|Y] \coloneqq \E_{Y=y} \H[ X | Y = y ]$$
where $(X|Y=y)$ is $X$ conditioned to the event $Y=y$ (again, we omit parentheses when there is no possibility of ambiguity).  The chain rule asserts that $\H[X|Y]$ can also be expressed by the formula
$$ \H[X|Y] = \H[X,Y] - \H[Y].$$
The \emph{mutual information} $\I(X:Y)$ between two random variables is given by the formula
$$ \I(X:Y) = \H[X] - \H[X|Y] = \H[Y] - \H(Y|X) = \H[X] + \H[Y] - \H[X,Y].$$
As is well known, $\I(X:Y)$ is non-negative, and vanishes precisely when $X,Y$ are independent.  Equivalently, one has the subadditivity property
$$ \H[X,Y] \leq \H[X] + \H[Y]$$
with equality precisely when $X,Y$ are independent.

We also define the conditional mutual information
$$ \I[X:Y|Z] \coloneqq \E_{Z=z} \I[(X|Z=z):(Y|Z=z)].$$
Clearly, $\I[X:Y|Z]$ is non-negative, and vanishes precisely when $X,Y$ are independent conditionally on $Z$.
From the chain rule we have
\begin{equation}\label{i-split}
\begin{split}
\I[X:Y|Z] &= \H[X|Z] - \H[X|Y,Z] \\
&= \H[Y|Z] - \H[Y|X,Z] \\
&= \H[X|Z] + \H[Y|Z] - \H[X,Y|Z].    
\end{split}
\end{equation}

\subsection{Acknowledgments}

The author was supported by the James and Carol Collins Chair, the Mathematical Analysis \& Application Research Fund, and by NSF grants DMS-2347850, and is particularly grateful to recent donors to the Research Fund.  He particularly thanks his coauthors Bogdan Georgiev, Javier G\'omez-Serrano, and Adam Zsolt Wagner for the highly productive and enjoyable collaboration \cite{gdm}, and for generously sharing the outputs of that collaboration for the purposes of writing the current paper.

While some of the results proven here were suggested by the outcome of AI-assisted experiments, the arguments in this paper are completely human-generated.

\section{A review of entropy sum set estimates}

Given two $\Q$-valued random variables $X,Y$ (not necessarily independent, or even defined on the same probability spaces), we define the \emph{entropic Ruzsa distance} $d[X;Y]$ between them by the formula
$$ d[X;Y] \coloneqq \H[X'-Y'] - \frac{1}{2} \H[X] - \frac{1}{2} \H[Y]$$
where $X',Y'$ are independent copies of $X,Y$ respectively.  This distance is non-negative, symmetric, and obeys the Ruzsa triangle inequality
\begin{equation}\label{ruzsa}
d[X;Z] \leq d[X;Y] + d[Y;Z]
\end{equation}
and the entropy bound
\begin{equation}\label{entropy-compare}
|\H[X]-\H[Y]| \leq 2d[X;Y];
\end{equation}
see \cite[Appendix A]{ggmt}. However we caution that this distance is not a true metric because $d[X;X]$ is non-zero in general.  We also observe the \emph{improved Ruzsa triangle inequality}
\begin{equation}\label{ruzsa-improv}
\H[X-Z] \leq \H[X-Y] + \H[Y-Z] - \H[Y]
\end{equation}
whenever $X,Y,Z$ are $\Q$-valued random variables with $Y$ independent of $X-Z$ (but $X,Z$ not required to be independent of each other); see \cite[Lemma 1.1]{gmt}.  We also recall the inequality
\begin{equation}\label{ruzsa-diff}
d[X;-Y] \leq 3d[X;Y],
\end{equation}
proven in \cite[Theorem 1.10]{tao-entropy}, as well as the Kaimanovich--Vershik--Madiman inequality
\begin{equation}\label{kvm}
\H[X+Y+Z] \leq \H[X+Y] + \H[X+Z] - \H[X]
\end{equation}
for any independent $\Q$-valued random variables $X,Y,Z$; see \cite[Lemma A.1]{ggmt}.

Now we study expressions of the form $\H[X-aX'] - \H[X]$ when $X'$ is an independent copy of $X$, where we see the role of rational complexity emerge.

\begin{proposition}\label{dilate}  Let $X$ be an $\Q$-valued random variable, and for any $a \in \Q$, let $g(a)$ denote the quantity $g(a) \coloneqq \H[X-aX'] - \H[X]$, where $X'$ is an independent copy of $X$.
\begin{itemize}
    \item[(i)] $g(0) = 0$ and $g(1)=d[X;X]$.
    \item[(ii)]  For any $a \in \Q$, we have $g(-a) \leq 3g(a)$, and if $a$ is non-zero, $g(a^{-1}) = g(a)$.
    \item[(iii)]  For any $k \geq 1$ and $a_1,\dots,a_k \in \Q$, we have $g(a_1 \dots a_k) \leq g(a_1)+ \dots + g(a_k)$ and $g(a_1 + \dots + a_k) \leq g(a_1)+\dots+g(a_k)+(k-1) g(1)$.
    \item[(iv)]  If $a$ is a non-zero integer, then $g(a) \leq (4 + 10 \lfloor \log_2 |a| \rfloor) g(1)$.
    \item[(v)] If $R$ is a finite set of slopes containing $0,1,\infty$ of cardinality $k+3$, then
    $$ g(a) \ll D(R;a)^{k+1} \max_{r \in R \backslash \{\infty\}} g(r).$$
\end{itemize}
\end{proposition}

\begin{proof}  The claim (i) is clear.  The first part of (ii) follows from \eqref{ruzsa-diff} for non-zero $a$ (noting that the case of $a=0$ follows from (i)), and the second part follows by observing that $-a (X - a^{-1} X')$ has the same distribution as $X - aX'$.  For (iii), it suffices by induction to treat the $k=2$ case.  For the first part of (iii), we observe from \eqref{ruzsa-improv} that
$$ \H[X - a_1 a_2 X'] \leq \H[X - a_1 X'] + \H[a_1X' - a_1a_2X''] - \H[a_1X']$$
when $X',X''$ are independent copies of $X$, which gives $g(a_1 a_2) \leq g(a_1)+g(a_2)$ after some changes of variable.  Similarly, from \eqref{ruzsa-improv} one has
$$ \H[X - a_1 X' - a_2X'] \leq \H[X - a_1 X' - X''] + \H[X'' - a_2X'] - \H[X'']$$
while from \eqref{kvm} one has
$$ \H[X - a_1 X' - X'']  \leq \H[X - a_1 X'] + \H[X-X''] - \H[X]$$
and the second part of (ii) follows by substituting the definition of $g$ and $d[X;X]$.

The claim (iv) follows from \cite[Lemma A.3(ii)]{ggmt-2}, but can also be derived from the previous portions of this proposition. (The constants in the bounds can be improved slightly if desired.)

Now we prove (v).  If we write $M \coloneqq \max_{r \in R \backslash \{\infty\}} g(r)$, then of course $g(r) \leq M$ for all $r \in R \backslash \{\infty\}$, in particular $g(1) \leq M$.  Hence by (iii), (iv) one has $g(m) \ll D(R;a) M$ whenever $m$ is a monomial in $R \backslash \{\infty\}$ of degree at most $D(R;a)$ and coefficient at most $2^{D(R;a)}$.  By (iii) again, we then have $g(p) \ll D(R;a)^{k+1} M$ whenever $p$ is a polynomial in the elements of $R \backslash \{\infty\}$ of degree at most $D(R;a)$ and coefficients at most $2^{D(R;a)}$.  Since $a$ is the ratio of two such polynomial expressions, the claim now follows from (ii), (iii). 
\end{proof}

We record the following form of the entropic Balog--Szemer\'edi--Gowers lemma.

\begin{lemma}[Entropic BSG]\label{entropic-bsg} Let $X,Y$ be $\Q$-valued random variables, and let $Z$ be a further random variable.  Let $(X',Y,Z)$ be a copy of $(X,Y,Z)$ with $X, X'$ conditionally independent relative to $(Y,Z)$.  Then
$$ \H[X-X'|Y,Z] \leq 2 \H[X+Y|Z] +\H[Y|Z] + 2\H[X|Z] - 2\H[X,Y|Z].$$
Also,
\begin{align*}
    &\E_{(X+Y,Z)=(c,z)} d[X|(X +Y,Z)=(c,z); Y|(X +Y,Z) = (c,z)] \\
    &\quad \leq 3\I[X : Y|Z]+2\H[X +Y|Z]-\H[X|Z]- \H[Y|Z].
\end{align*}
\end{lemma}

\begin{proof} By averaging in $Z$, it suffices to prove these two inequalities in the unconditional case when $Z$ is not present.  But then this follows from (the proof of) \cite[Lemma 3.3]{tao-entropy} and \cite[Lemma A.2]{ggmt} respectively.
\end{proof}

\section{A reduction to the independent case}

We consider the following variant of $\SD(R;s)$, in which the random variables $X,Y$ are now assumed to be independent, which is easier to control by the ``entropic Pl\"unnecke--Ruzsa calculus'' of the previous section than the original sum-difference constants.

\begin{definition}[Independent sum-difference constant]\label{fdef}  Let $R$ be a finite set of slopes containing $\{0,1,\infty\}$, and let $s$ be a further slope not in $R$.  Then $f(R;s)$ is defined to be the least constant for which the following statement holds (or $+\infty$ if no such constant exists): whenever $X,Y$ are independent $\Q$-valued random variables, and $M \geq 0$, $K \geq 1$ are such that
\begin{equation}\label{mk}
 M - \log K \leq \H[\pi_r(X,Y)] \leq M + \log K
 \end{equation}
for all $r \in R$, then
$$\H[\pi_{s}(X,Y)] \leq M + f(R;s) \log K.$$
\end{definition}

If one sets $M \coloneqq \H[X]$ and $\log K \coloneqq \sum_{r \in R} |\H[\pi_r(X,Y)] - \H[X]|$, we see in particular that
\begin{equation}\label{fdef-alt}
\H[\pi_{s}(X,Y)] \leq \H[X] + f(R;s) \sum_{r \in R} |\H[\pi_r(X,Y)] - \H[X]|
\end{equation}
whenever $X,Y$ are independent.

If we let $X,Y$  be independent and uniform on $\{1,\dots,K\}$ for a large integer $K$, then $\H[\pi_r(X,Y)] = \log K + O(1)$ for any given slope $r$ (with implied constants depending on $r$).  Sending $K \to \infty$ and setting $M=0$, we conclude that 
\begin{equation}\label{frs-lower}
f(R;s) \geq 1.
\end{equation}
The main result of this section is to relate $f(R;s)$ to $\SD(R;s)$ as follows.

\begin{theorem}[Reduction to independent case]\label{indep} Let $R$ be a finite set of slopes with $0, \infty \in S$, and let $s$ be a further slope not in $S$.  Then
$$ \SD(R;s) \leq 2 - \frac{c}{f(R;s)}$$
for some absolute constant $c>0$.
\end{theorem}

\subsection{First step: obtaining a good configuration}  
We now prove the theorem.  
In principle this is ``just'' an application of the (entropic) Balog--Szemer\'edi--Gowers theorem, Lemma \ref{entropic-bsg}, to replace the coupled random variables $X,Y$ by independent versions, but (similarly to \cite{ bourgain}, though now in the entropic framework) one has to take some care not to obtain control on various doubling constants and projections when doing so.  As such, the proof shall be somewhat lengthy and proceed in stages.

We may assume without loss of generality that $f(R;s)$ is finite, since the claim is trivial otherwise. (The finiteness of $f(R;s)$ will be demonstrated in the next section, but we will not need that result for the current argument.)
Suppose that $\SD(R;s) > 2 - \delta$ for some small $\delta > 0$.  We will show that 
\begin{equation}\label{fF-targ}
f(R;s) \delta \gg 1, 
\end{equation}
which will imply the claim.

Our argument shall involve the following key definitions.

\begin{definition}[Configurations]\label{config-def}  A \emph{configuration} is a quadruple $(M,X,Y,Z)$, where $M > 0$, and $X,Y,Z$ are random variables (not necessarily independent) with $X,Y$ $\Q$-valued.
\begin{itemize}
    \item[(i)] A configuration is \emph{good} if it obeys the bounds
\begin{equation}\label{pir-1}
\H[\pi_r(X,Y)|Z] \leq M
\end{equation}
for all $r \in R$, and also
\begin{equation}\label{hpis}
 \H[\pi_{s}(X,Y)|Z] \geq (2-O(\delta))M.
 \end{equation}
    \item[(ii)] A configuration $(M,X,Y,Z)$ has an \emph{invertible projection} if, for each $z$ in the essential range of $Z$, the map $\pi_s$ is injective
on the essential range of $((X,Y)|Z=z)$.
    \item[(iii)] A configuration $(M,X,Y,Z)$ has \emph{controlled $X$-doubling} if one has
\begin{equation}\label{x-cont}
\E_{Z=z} d[X|Z=z; X|Z=z] \ll \delta M.
\end{equation}
    \item[(iv)] A configuration $(M,X,Y,Z)$ has \emph{controlled $Y$-doubling} if one has
\begin{equation}\label{y-cont}
\E_{Z=z} d[Y|Z=z; Y|Z=z] \ll \delta M.
\end{equation}
\end{itemize}
\end{definition}

By definition of $\SD(R;s)$, one can find $\Q$-valued random variables $X,Y$ and some $M > 0$ such that
$$ \H[\pi_r(X,Y)] \leq M$$
for all $r \in R$, and
$$ \H[\pi_{s}(X,Y)] \geq (2-\delta) M.$$
Thus there exists at least one good configuration $(M,X,Y,Z)$ (with $Z$ trivial in this case).  However, this configuration is not guaranteed to have an invertible projection or have controlled $X$-doubling or $Y$-doubling.

In the next few steps of the arguments, we shall locate configurations with more of the properties listed in Definition \ref{config-def}.  For now, we record some basic properties of good configurations.

\begin{lemma}[Properties of good configurations]\label{good}  Let $(M,X,Y,Z)$ be a good configuration.  Then
\begin{equation}\label{admis-1}
\H[\pi_r(X,Y)|Z] = (1 + O(\delta)) M
\end{equation}
for all $r \in R$ and
\begin{equation}\label{admis-2}
\H[X,Y|Z], \H[\pi_{s}(X,Y)|Z] = (2 + O(\delta)) M.
\end{equation}
In particular we have
\begin{equation}\label{admis-3}
\H[X|Z], \H[Y|Z] = (1 + O(\delta)) M
\end{equation}
and hence by \eqref{i-split}
\begin{equation}\label{admis-4}
\I[X:Y|Z] \ll \delta M.
\end{equation}
\end{lemma}

\begin{proof}
For any two distinct slopes $r,r' \in \R$, $\pi_r(X,Y)$ and $\pi_{r'}(X,Y)$ uniquely determine $(X,Y)$, hence
so in particular
$$ \H[\pi_{s}(X,Y)|Z] \leq \H[X,Y|Z] \leq \H[\pi_r(X,Y)|Z] + \H[\pi_{r'}(X,Y)|Z] \leq 2M.$$
Comparing this with \eqref{hpis}, we obtain the claim.
\end{proof}

\subsection{Second step: ensuring an invertible projection}

It will be convenient to work with configurations $(M,X,Y,Z)$ with an invertible projection.  We record some basic properties of such configurations:

\begin{lemma}[Consequences of an invertible projection]\label{invert-consequence}  If $(M,X,Y,Z)$ is a configuration with an invertible projection, then
\begin{equation}\label{hpisxyz}
 \H[\pi_s(X,Y)|Z] = \H[X,Y|Z].
\end{equation}
Also, $(M,X,Y,(W,Z))$ has an invertible projection for any further random variable $W$.
\end{lemma}

\begin{proof} For any $z$ in the essential range of $Z$, the invertibility of $\pi_s$ on $((X,Y)|Z=z)$ implies that
$$ \H[\pi_s(X,Y)|Z=z] = \H[X,Y|Z=z].$$
Averaging over $z$, we obtain the first claim.  The second claim is trivial, since the essential range of $((X,Y)|(W,Z)=(w,z))$ is contained in that of $((X,Y)|Z=z)$ for any $(w,z)$ in the essential range of $(W,Z)$.
\end{proof}

The main result of this step is

\begin{proposition}[Ensuring invertible projections]\label{doubt}  There exists a good configuration $(M,X,Y,Z)$ with an invertible projection.
\end{proposition}

\begin{proof}  
This type of reduction is easy in the combinatorial setting \eqref{pise}, as one simply removes all but one element from $E$ in each fiber of $\pi_s$.  To do the analogous operation in the entropic setting we will need to perform this removal operation randomly, as follows.

By the previous discussion, we can find a good configuration $(M,X,Y,0)$ with the ``$Z$'' variable trivial.  Let $C$ denote the essential range of $\pi_s(X,Y)$.  For each $c \in C$, let $\Gamma_c$ denote the essential range of the conditioned random variable $((X,Y)|\pi_s(X,Y)=c)$.  This is a finite non-empty subset of $\{(a,b) \in \Q \times \Q: \pi_{s}(a,b) = c \}$.  We then define a random function $\phi \colon C \to \Q \times \Q$ by declaring $\phi(c)$ for $c \in C$ to be an element of $\Gamma_c$ chosen with the distribution of $((X,Y)|\pi_{s}(X,Y)=c)$, with the random variables $\phi(c)$ being jointly independent in $c$, and also independent of $(X,Y)$.  We then define a ``shuffled'' version $(X',Y')$ of $(X,Y)$ by the formula
$$ (X',Y') \coloneqq \phi( \pi_{s}(X,Y) ).$$
Observe that on the event $\pi_{s}(X,Y)=c$, the random variables $(X,Y)$ and $(X',Y')$ have identical distribution by construction, and hence $(X',Y')$ also unconditionally has the same distribution as $(X,Y)$.  In particular, for any $r \in R$ we have
$$ \H[\pi_r(X',Y')|\phi] \leq \H[\pi_r(X',Y')] = \H[\pi_r(X,Y)] \leq M.$$
Also, it is clear from construction that $\pi_s$ is injective on the essential range of $((X',Y')|\phi=\phi_0)$ for every $\phi_0$ in the essential range of $\phi$, and $\pi_s(X',Y') = \pi_s(X,Y)$ almost surely. In particular,
$$ \H[\pi_s(X',Y')|\phi] = \H[\pi_s(X,Y)|\phi] = \H[\pi_s(X,Y)] \geq (2-\delta) M.$$

Thus $(M, X', Y', \phi)$ is a good configuration with an invertible projection, giving the claim.
\end{proof}

\subsection{Third step: controlling one doubling constant}

Now, we use the entropic Balog--Szemer\'edi--Gowers lemma to ensure controlled $X$-doubling.

\begin{proposition}[Ensuring controlled $X$-doubling]\label{xdub}  Let $(M,X,Y,Z)$ be a good configuration with an invertible projection.  Then there exists an additional random variable $W$ such that
$(M,X,Y,(W,Z))$ is a good configuration with an invertible projection and controlled $X$-doubling.
\end{proposition}

\begin{proof}
Let $(Y',X,Z)$ be a copy of $(Y,X,Z)$ with $Y, Y'$ conditionally independent relative to $(X,Z)$, then we also have
$$ \H[X|Z], \H[Y'|Z], \H[X+Y'|Z] = (1 + O(\delta)) M; \quad \H[X,Y'|Z] = (2+O(\delta)) M.$$
Let $(X',Y',Z)$ be a copy of $(X,Y',Z)$ with $X, X'$ conditionally independent relative to $(Y',Z)$.  By the first part of Lemma \ref{entropic-bsg} we have
$$ \H[X-X'|Y',Z] \leq (1 + O(\delta)) M.$$
We also have
$$ \H[X|Y',Z] = \H[X,Y'|Z] - \H[Y'|Z] = (1 + O(\delta)) M.$$
We conclude that
\begin{align*}
\E_{(Y',Z)=(y,z)} &d[X|(Y',Z)=(y,z); X|(Y',Z)=(y,z)]\\
&= \H[X-X'|Y'-Z] - \frac{1}{2} \H[X|Y',Z] - \frac{1}{2} \H[X|Y',Z] \\
&\ll \delta M.    
\end{align*}
Thus $(M,X,Y,(Y',Z))$ has controlled $X$-doubling.
For any $r \in R$, we have
\begin{equation}\label{prx}
 \H[\pi_r(X,Y)|Y',Z] \leq \H[\pi_r(X,Y)|Z] \leq M.
\end{equation}
Also, we have
\begin{align*}
    \H[X,Y|Y',Z] &= \H[X,Y,Y'|Z] - \H[Y'|Z] \\
    &= \H[X,Y|Z] + \H[X,Y'|Z] - \H[X|Z] - \H[Y'|Z] \\
    &= (2 + O(\delta)) M.
\end{align*}
By Lemma \ref{invert-consequence}, $(M,X,Y,(Y',Z))$ has an invertible projection and
$$ \H[\pi_s(X,Y)|Y',Z] = \H[X,Y|Y',Z] = (2 + O(\delta)) M.$$
We conclude that $(M,X,Y,(Y',Z))$ is a good configuration.  The claim follows.
\end{proof}

\subsection{Fourth step: controlling a second doubling constant}

By swapping the roles of $X$ and $Y$ in Proposition \ref{xdub}, we have

\begin{proposition}[Ensuring controlled $Y$-doubling]\label{ydub}  Let $(M,X,Y,Z)$ be a good configuration with an invertible projection.  Then there exists an additional random variable $W$ such that
$(M,X,Y,(W,Z))$ is a good configuration with an invertible projection and controlled $Y$-doubling.
\end{proposition}

We would like to combine Proposition \ref{xdub} with Proposition \ref{ydub} to obtain a good configuration with invertible projections that simultaneously has controlled $X$-doubling and controlled $Y$-doubling.  However, there is a technical issue: the passage of 
$(M,X,Y,Z)$ to $(M,X,Y,(W,Z))$ in Proposition \ref{ydub} could potentially destroy the controlled $X$-doubling property.  Fortunately, this does not occur, thanks to the following lemma:

\begin{lemma}\label{sense}  Suppose that $(M,X,Y,Z)$ and $(M,X,Y,(W,Z))$ are both good configurations.  If $(M,X,Y,Z)$ has controlled $X$-doubling, then so does $(M,X,Y,(W,Z))$.
\end{lemma}

\begin{proof}
Let $(\tilde X,Z)$ be a copy of $(X,Z)$ with $\tilde X$, $(X,W)$ conditionally independent over $Z$, which also implies that $X, \tilde X$ are conditionally independent over $(W,Z)$.  By the Ruzsa triangle inequality \eqref{ruzsa}, we have
\begin{align*}
\E_{(W,Z)=(w,z)} &d[X|(W,Z)=(w,z); X|(W,Z)=(w,z)] \\
&\leq
2 \E_{(W,Z)=(w,z)} d[X|(W,Z)=(w,z); \tilde X|(W,Z)=(w,z)]  \\
&= 2 \left(\H[X-\tilde X|W,Z] - \frac{1}{2}\H[X|W,Z] - \frac{1}{2} \H[\tilde X|W,Z]\right) \\
&\leq 2 \left(\H[X-\tilde X|Z] - \frac{1}{2}\H[X|W,Z] - \frac{1}{2} \H[X|Z]\right) \\
&\leq 2 \left(\H[X-\tilde X|Z] - \frac{1}{2}\H[X|Z] - \frac{1}{2} \H[X|Z]\right) + O(\delta M) \\
&= 2 \E_{Z=z} d[X|Z=z; X|Z=z] + O(\delta M) \\
&\ll \delta M
\end{align*}
thanks to \eqref{admis-3} and \eqref{x-cont}. Thus $(M,X,Y,(W,Z))$ has controlled $X$-doubling as required.
\end{proof}

Combining Proposition \ref{xdub} with Proposition \ref{ydub} and Lemma \ref{sense}, we now conclude

\begin{corollary}[Simultaneous controlled doubling]\label{xydub}  There exists a good configuration $(M,X,Y,Z)$ with an invertible projection that has both controlled $X$-doubling and controlled $Y$-doubling.
\end{corollary}

\subsection{Fifth step: controlling projections of conditionally independent variables}

The configuration $(M,X,Y,Z)$ provided by Corollary \ref{xydub} obeys many further estimates:

\begin{proposition}  Let $(M,X,Y,Z)$ be as in Corollary \ref{xydub}. Let $r \in R \backslash \{0,\infty\}$. Then we have the bounds
\begin{align}
\E_{(X+rY,Z)=(c,z)} d[X|(X+rY,Z)=(c,z); rY|(X+rY,Z)=(c,z)] &\ll \delta M \label{s1}\\
\E_{(X+rY,Z)=(c,z)} d[X|(X+rY,Z)=(c,z); X|Z=z] &\ll \delta M\label{s2}\\
\E_{(X+rY,Z)=(c,z)} d[rY|(X+rY,Z)=(c,z); rY|Z=z] &\ll \delta M.\label{s3}
\end{align}
In particular, from the Ruzsa triangle inequality \eqref{ruzsa} one has
$$\E_{(X+rY,Z)=(c,z)} d[X|Z=z; rY|Z=z] \ll \delta M$$
and hence
\begin{equation}\label{s4}
\E_{Z=z} d[X|Z=z; rY|Z=z] \ll \delta M.
\end{equation}
\end{proposition}

\begin{proof}  We begin with \eqref{s1}.  By the second part of Lemma \ref{entropic-bsg}, the left-hand side is bounded by
$$ 3 \I[X:rY|Z] + 2\H[X+rY|Z] - \H[X|Z] - \H[Y|Z].$$
The claim now follows from Lemma \ref{good}.

Now we establish \eqref{s2}, which will be proven by a variant of the argument used to establish Lemma \ref{sense}.  Let $(X',Z)$ be a copy of $(X,Z)$ with $X', (X,Y)$ independent conditionally on $Z$, which implies that $X'$ is independent of $X+rY$ conditionally on $Z$.  We compute
\begin{align*}
\E_{(X+rY,Z)=(c,z)}& d[X|(X+rY,Z)=(c,z); X|Z=z]\\ 
&= \E_{(X+rY,Z)=(c,z)} d[X|(X+rY,Z)=(c,z); X'|(X+rY,Z)=(c,z)] \\
&= \H[X-X'|X+rY,Z] - \frac{1}{2} \H[X|X+rY,Z] - \frac{1}{2} \H[X'|X+rY,Z] \\
&\leq \H[X-X'|Z] - \frac{1}{2} \left(\H[X,X+rY|Z]-\H[X+rY|Z]\right) - \frac{1}{2} \H[X'|Z] \\
&\leq \H[X-X'|Z] - \frac{1}{2} \H[X|Z] - \frac{1}{2} \H[X'|Z] + O(\delta M)\\
&= \E_{Z''=z''} d[X|Z=z; X|Z=z] + O(\delta M)\\
&\ll \delta M    
\end{align*}
as claimed thanks to Lemma \ref{good} and \eqref{x-cont}.

To prove \eqref{s3}, we similarly let $(Y',Z)$ be a copy of $(Y,Z)$ with $Y', (X,Y)$ independent conditionally on $Z$, which implies that $rY'$ is independent of $X+rY$ conditionally on $Z$.  Then
\begin{align*}
 \E_{(X+rY,Z)=(c,z)}& d[rY|(X+rY,Z)=(c,z); rY|Z=z] \\
&= \E_{(X+rY,Z)=(c,z)} d[rY|(X+rY,Z)=(c,z); rY'|(X+rY,Z)=(c,z)] \\
&= \H[rY-rY'|X+rY,Z] - \frac{1}{2} \H[rY|X+rY,Z] - \frac{1}{2} \H[rY'|X+rY,Z] \\
&\leq \H[rY-rY'|Z] - \frac{1}{2} \left(\H[rY,X+rY|Z]-\H[X+rY|Z]\right) - \frac{1}{2} \H[rY'|Z] \\
&\leq \H[Y-Y'|Z] - \frac{1}{2} \H[Y|Z] - \frac{1}{2} \H[Y'|Z] + O(\delta M)\\
&= \E_{Z=z} d[Y|Z=z; Y|Z=z] + O(\delta M)\\
&\ll \delta M    
\end{align*}
thanks to Lemma \ref{good} and \eqref{y-cont}.
\end{proof}

\subsection{Final step: using the quantity $f(R;s)$}

Let $(M,X,Y,Z)$ be as in Corollary \ref{xydub}.  
If $(Y',Z)$ is a copy of $(Y,Z)$ with $Y'$ independent of $(X,Y)$ conditionally on $Z$, we conclude from \eqref{s4} and \eqref{entropy-compare} that
$$ \E_{Z=z} |\H[X+rY'|Z=z] - \H[X|Z=z]| \ll \delta M$$
for all $r \in R \backslash \{0,\infty\}$.  
Applying \eqref{fdef-alt} and the triangle inequality, we conclude that
$$
 \E_{Z=z} \H[X+sY'|Z=z] \leq \E_{Z=z} \H[X|Z=z] + O(f(R;s) \delta M).
$$
From Lemma \ref{good} we have
$$ \E_{Z=z} \H[X|Z=z] = \H[X|Z] = (1+O(\delta)) M$$
and thus (by \eqref{frs-lower})
$$ \H[X+sY'|Z]  =  \E_{Z=z} \H[X+sY'|Z=z] \leq (1 + O(f(R;s)\delta)) M.$$
Comparing this with \eqref{admis-2} and using \eqref{frs-lower}, we obtain \eqref{fF-targ} as desired.

\section{Proof of Theorem \ref{main-2}}\label{main-sec}

We now prove Theorem \ref{main-2}.  We begin with part (i).  Write $s=a/b$, where $b \geq 1$ and $a \neq 0,b$ is coprime to $b$, then $D(R;s) \asymp \log(2+|a|+|b|)$.  Let $N$ be the natural number $\max( 2, |a|/10, b/10)$.  If one takes $X,Y$ to be independent copies of the uniform distribution on $\{1,\dots,N\}$, we can compute that $\H[X], \H[Y], \H[X+Y] = \log N + O(1)$, and $\H[X+\frac{a}{b}Y] = 2 \log N + O(1)$ (noting that $\pi_{a/b}$ is injective on $\{1,\dots,N\} \times \{1,\dots,N\}$ if $N>2$).  This implies that
$$\SD\left(\{0,1,\infty\};\frac{a}{b}\right) \geq 2 - O\left(\frac{1}{\log N}\right)$$
which gives the required lower bound.  For the matching upper bound, it suffices by Theorem \ref{indep} to show that $f(\{0,1,\infty\},a/b) \ll \log(2+|a|+|b|)$.  That is to say, given independent random variables $X,Y$ obeying bounds of the form
$$ \H[X], \H[Y], \H[X+Y] = M + O(\log K)$$
for some $M>0$ and $K \geq 2$, we need to show that
$$ \H\left[X + \frac{a}{b} Y\right] = M + O (\log(2+|a|+|b|) \log K).$$
From the Ruzsa triangle inequality one has
$$ d[X;X] \ll \log K.$$
Applying Proposition \ref{dilate}(iv), we have
$$ \H[X + aY] = M + O(\log(2+|a|) \log K)$$
and
$$ \H[X - bX] = M + O(\log(2+|b|) \log K)$$
and hence by the Ruzsa triangle inequality \eqref{ruzsa} (or \eqref{ruzsa-improv}) we have
$$ \H[bX + aY] = M + O(\log(2+|a|+|b|) \log K),$$
giving the claim.

We now prove part (ii).  We allow implied constants to depend on $k$, and abbreviate $D(R;s)$ as $D$.  We first show the lower bound
$$ \SD(R; s) \geq 2 - O\left(\frac{\log(2+D)}{D} \right)$$
whenever $s$ has complexity at least $D$.  Clearly we may assume that $D$ is larger than any constant depending on $k$.

For any $N \geq 1$, let $P_N \subset \Q$ be the set of all polynomial combinations of slopes in $R \backslash \{0,1,\infty\}$ of degree at most $N$ with coefficients between $-2^{N-1}$ and $2^{N-1}$.  By crude counting we have
$$ \log |P_{D}| \leq \log\left( (2^{D}+1)^{(D+1)^{k}} \right) \ll D^{k+1}.$$
Thus by the pigeonhole principle, one can find $1 \leq N < D$ such that
$$ \log |P_{N+1}| \leq \left(1 + O\left(\frac{\log D}{D} \right)\right) \log |P_N|.$$
If we then let $X,Y$ be drawn uniformly at random from $P_N$, then for any $r \in R$, $\pi_r(X,Y)$ takes values in $P_{N+1}$, hence from Jensen's inequality we have
$$ \H[\pi_r(X,Y)] \leq \log |P_{N+1}| \leq \left(1 + O\left(\frac{\log D}{D} \right)\right) \log |P_N|.$$
On the other hand, if $\pi_s$ is not injective on the range of $(X,Y)$, so that
$$ p_1 + s p_2 = p_3 + s p_4$$
for some polynomial combinations $p_1,p_2,p_3,p_4$ of slopes in $R \backslash \{\infty\}$ of degree at most $D-1$ and coefficients of magnitude at most $2^{D-2}$ with $p_2 \neq p_4$, then we may rearrange as
$$ s = \frac{p_3 - p_1}{p_2 - p_4}$$
and conclude that $s$ has rational complexity at most $D-1$, a contradiction.  Thus $\pi_s$ is injective and
$$ \H[\pi_s(X,Y)] = \H[X,Y] = 2 \log |P_N|.$$
Inserting these bounds into \eqref{sdef}, we obtain the claim.

To show the upper bound, it suffices by Theorem \ref{indep} to show that when $s$ has complexity $D$, that
$$  f(R;s) \ll D^{k+1}.$$
Thus, given independent random variables $X,Y$ obeying bounds of the form
\begin{equation}\label{hx}
\H[X], \H[Y], \H[X+Y], \H[X+rY] = M + O(\log K)
\end{equation}
for some $M>0$ and $K \geq 2$, and all $r \in R \backslash \{0,1,\infty\}$, it will suffice to show that
\begin{equation}\label{dxsy}
d[-X; sY] \ll D^{k+1} \log K.
\end{equation}
Note from \eqref{ruzsa}, \eqref{ruzsa-diff} that 
\begin{equation}\label{dxy}
d[X;Y], d[Y;Y] \ll \log K.
\end{equation}

For $a \in \Q$, define $g(a) \coloneqq \H[Y-aY'] - \H[Y]$, where $Y'$ is an independent copy of $Y$. From \eqref{hx}, \eqref{ruzsa}, \eqref{ruzsa-diff} we have $g(r) \ll \log K$ for $r \in R \backslash \{\infty\}$.  By Proposition \ref{dilate}(iv), we conclude that $g(s) \ll D^{k+1} \log K$, thus $d[Y; sY] \ll D^{k+1} \log K$.  Combining with \eqref{dxy} and \eqref{ruzsa}, \eqref{ruzsa-diff}, we obtain \eqref{dxsy} as required.

\section{A continuous limit}\label{cts-sec}

In this section we establish \Cref{cts-limit}.  Let $\alpha$, $a$, $b$, $f$, $f_0, f_1, f_\infty$ be as in that theorem.  We allow implied constants in the asymptotic notation to depend on $\alpha, f$. We will establish the lower bound
\begin{equation}\label{sd-targ}
 \SD\left(\{0,1,\infty\},\frac{a}{b}\right) \geq 2 - \frac{h_2(f_\alpha) - \max(h(f_0), h(f_\infty), h(f_1))+o(1)}{\log b}.
 \end{equation}
Comparing this with Theorem \ref{main-2}, we conclude that the expression $h_2(f_\alpha) - \max(h(f_0), h(f_\infty), h(f_1))$ is bounded from below, hence $c_\alpha > 0$; by taking an arbitrary test function for $f$ we also see that $c_\alpha < \infty$.  Taking $h_2(f_\alpha) - \max(h(f_0), h(f_\infty), h(f_1))$ arbitrarily close to $c_\alpha$, we obtain the desired claim \eqref{up}.

It remains to establish \eqref{sd-targ}.
By the Poisson summation formula, the rapid decrease of the Fourier transform of the smooth compactly supported $f$, and the mass one hypothesis, we have
\begin{equation}\label{bib}
 \frac{1}{b^2} \sum_{n,m \in \Z} f\left(\frac{n}{b}, \frac{m}{b}\right) = 1 + \eps
 \end{equation}
for some $\eps = O(1/b)$; in fact one can get much better decay than this, but for our purposes any decay faster than $1/\log b$ will suffice.  We then take $X,Y$ to be supported on the grid $\Z^2$ with probability distribution
$$ \P((X,Y) = (n,m)) \coloneqq \frac{1}{(1+\eps) b^2} \left(\frac{n}{b}, \frac{m}{b}\right)$$
thus $(X,Y)$ takes values in a ball of radius $O(b)$  and
\begin{equation}\label{pnx}
\P((X,Y) = (n,m)) = \frac{1}{b^2} f\left(\frac{n}{b}, \frac{m}{b}\right) + O\left( \frac{1}{b^3} \right).
\end{equation}
By Bezout's theorem, any integer can be uniquely written in the form $bn+am$ for some $0 \leq m < b$, and any other representation of the form $bn'+am'$ takes the form $b(n-aj) + a(m+bj)$.  Thus
$$ \P(bX+aY = bn+am) = \sum_{j \in\Z} \P( (X,Y) = (n-aj, m+bj) ).$$
There are only $O(1)$ values of $j$ for which this sum is non-zero, so from \eqref{pnx} one has
$$ \P(bX+aY = bn+m) = \frac{1}{b^2} f_{(a/b)}\left(\frac{n}{b}, \frac{m}{b}\right) + O\left( \frac{1}{b^3} \right).$$
where $f_{(a,b)}: \R \times [0,1) \to \R^+$ is the function
$$f_{(a/b)}(x,y) \coloneqq \sum_{j \in \Z} f\left(x-\frac{a}{b} j, y+j\right).$$
Applying the entropy function $\h$, we conclude that
$$ \h(\P(bX+aY = bn+am)) = \frac{2 \log b}{b^2} f_{(a/b)}\left(\frac{n}{b}, \frac{m}{b}\right) +
\frac{1}{b^2} \h\left(f_{(a/b)}\left(\frac{n}{b}, \frac{m}{b}\right)\right) + O\left( \frac{\log b}{b^3} \right).$$
From \eqref{bib} one has
$$ \frac{1}{b^2} \sum_{n \in\Z; 0 \leq m < b} f_{(a/b)}\left(\frac{n}{b}, \frac{m}{b}\right) = 1 + \eps = 1 + O\left(\frac{1}{b} \right)$$
and from (uniform) Riemann integrability of the $f_{(a/b)}$ one has
$$ \frac{1}{b^2} \sum_{n \in\Z; 0 \leq m < b} \h\left(f_{(a/b)}\left(\frac{n}{b}, \frac{m}{b}\right)\right) = h_2(f_{(a/b)}) + o(1).$$
Finally, from dominated convergence one has
$$ h_2(f_{(a/b)}) = h_2(f_{(\alpha)})+o(1)$$
so we conclude that
$$ \H\left[\pi_{a/b}(X,Y)\right] = \H[bX+aY]  = 2 \log b + h_2(f_{(\alpha)}) + o(1).$$
In a similar vein, from another application of Poisson summation and \eqref{f0-def} we see that
$$ \P(X=n) = \frac{1}{b} f_0\left(\frac{n}{b}\right) + O\left(\frac{1}{b^2}\right)$$
for any integer $n$, hence
$$ \h(\P(X=n)) = \frac{\log b}{b} f_0\left(\frac{n}{b}\right) + \frac{1}{b} \h\left(f_0 \left(\frac{n}{b}\right)\right) + O\left(\frac{\log b}{b^2}\right).$$
Since $X = O(b)$, one can sum using Riemann integrability and \eqref{bib} to conclude that
$$ \H[\pi_0(X,Y)] = \H[X]= \log b + h(f_0) + o(1).$$
Similar arguments give
$$ \H[\pi_\infty(X,Y)] = \H[Y] = \log b + h(f_\infty) + o(1)$$
and
$$ \H[\pi_1(X,Y)] = \H[X+Y] = \log b + h(f_1) + o(1).$$
Comparing this with \eqref{sdef}, we obtain \eqref{sd-targ}.

\end{document}